\documentclass{icgg}

\usepackage{graphicx}     % if you want to include images
\usepackage{microtype}    % improves, among other things, line-breaks
\usepackage{dsfont,amsfonts,amsthm}
\usepackage{nicefrac}
\usepackage{natbib}
\allowdisplaybreaks

\newtheorem{definition}{Definition}
\newtheorem{remark}{Remark}
\newtheorem{theorem}{Theorem}
\newtheorem{example}{Example}

\DeclareMathAlphabet{\mathcal}{OMS}{cmsy}{m}{n}
\newcommand{\clifford}{\mathcal{C}\ell}
\newcommand{\corr}{\ensuremath{\mathrel{\widehat{=}}}} 

\title{NULL POLARITIES AS GENERATORS OF THE PROJECTIVE GROUP}

\author{Daniel KLAWITTER
 }
\affiliation{Technical University of Dresden, Institute of Geometry, Germany}

\hypersetup{
  pdfauthor = {Daniel Klawitter},
  pdftitle = {Null Polarities as Generators of the Projective Group},
  pdfsubject = {paper003},
  pdfkeywords = {Geometry, graphics, ISGG, ICGG, LaTeX, article
    formatting}
}

\begin{document}
% Set papernumber as received upon acceptance of your contribution.
\papernumber{003}

% We need a little trickery in order to get a one-column abstract in a
% two-column text:
\twocolumn[%
\begin{@twocolumnfalse}
  \maketitle
  % Here comes the abstract...
\begin{abstract}
It is well-known that the group of regular projective transformations of $\mathbb{P}^3(\mathbb{R})$ is isomorphic to the group of projective automorphisms of Klein's quadric $M_2^4\subset\mathbb{P}^5(\mathbb{R})$. We introduce the Clifford algebra $\clifford_{(3,3)}$ constructed over the quadratic space $\mathbb{R}^{(3,3)}$ and describe how points on Klein's quadric are embedded as null vectors, {\it i.e.},  grade-$1$ elements squaring to zero. Furthermore, we discuss how geometric entities from Klein's model can be transferred to this homogeneous Clifford algebra model. Automorphic collineations of Klein's quadric can be described by the action of the so called sandwich operator applied to vectors $\mathfrak{v}\in\bigwedge^1 V$. Vectors correspond to null polarities in $\mathbb{P}^3(\mathbb{R})$. We introduce a factorization algorithm. With the help of this algorithm we are able to factorize an arbitrary versor $\mathfrak{g}\in\clifford_{(3,3)}$ into a set of non-commuting vectors $\mathfrak{v}_i\in\bigwedge^1 V,\,i=1,\dots, k,\, 1\leq k\leq 6$ corresponding to null polarities with $\mathfrak{g}=\mathfrak{v}_1\dots\mathfrak{v}_k$. Thus, we present a method to factorize every collineation in $\mathbb{P}^5(\mathbb{R})$ that is induced by a projective transformation acting on $\mathbb{P}^3(\mathbb{R})$ into a set of at most six involutoric automorphic collineations of Klein's quadric corresponding to null polarities respectively skew-symmetric $4 \times 4$ matrices.
Moreover, we give an outlook for Lie's sphere geometry, {\it i.e.}, the homogeneous Clifford algebra model constructed with the quadratic form corresponding to Lie's quadric $L_1^{n+1}\subset\mathbb{P}^{n+2}(\mathbb{R})$.
\end{abstract}
  \keywords{Clifford algebra, line geometry, Klein's quadric, null polarity, factorization}
\end{@twocolumnfalse}]

% Starting from here, everything is, more or less, standard LaTeX...

\section{Introduction}
Klein's quadric denoted by $M_2^4\subset\mathbb{P}^5(\mathbb{R})$ and subspace intersections of Klein's quadric are well-known, cf \cite{potwal} or \cite{weiss:vor}. The group of regular projective transformations of $\mathbb{P}^3(\mathbb{R})$ is isomorphic to the group of automorphic collineations of Klein's quadric. Moreover, Cayley-Klein geometries can be represented by Clifford algebras, where the group of Cayley-Klein isometries is given by the Pin group of the corresponding Clifford algebra, see \cite{gunn:kinematics}. Therefore, we recall fundamental basics of Clifford algebras and recall the homogeneous Clifford algebra model corresponding to Klein's quadric. The relationsship between projective transformations and elements of the Pin group is given. The full construction can be found in \cite{klawitter_lingeo}.

\section{Geometric Algebra}
General introductions to Clifford algebras can be found for example in \cite{gallier:cliffordalgebrascliffordgroups,garling:cliffordalgebras,perwass} and \cite{porteous:cliffordalgebrasandtheclassicalgroups}.

\begin{definition}
Let $V$ be a real valued vector space of dimension $n$. Furthermore, let $b:V\mapsto \mathbb{R}$ be a quadratic form on $V$. The pair $(V,b)$ is called \emph{quadratic space}.
\end{definition}
\noindent
We denote the Matrix corresponding to $b$ by $\mathrm{B}_{ij}$ with $i\leq j,j\leq n$. Therefore $b(x_i,x_j)=\mathrm{B}_{ij}$ for some basis vectors $x_i$ and $x_j$.
\begin{definition}
The Clifford algebra is defined by the relations
\begin{equation}
\label{EQ1}
x_ix_j+x_jx_i=2\mathrm{B}_{ij},\quad 1\leq i,j\leq n.
\end{equation}
\end{definition}
\noindent
Usually, the algebra is denoted by $\clifford(V,b)$.
By Silvester's law of inertia we can always find a basis $\lbrace e_1,\dots,e_n \rbrace$ of $V$ such that $e_i^2$ is either $1,-1$ or $0$.
\begin{definition}
The number of basis vectors that square to $(1,-1,0)$ is called \emph{signature} $(p,q,r)$. If $r\neq 0$ we call the geometric algebra \emph{degenerated}. We denote a Clifford algebra by $\clifford_{(p,q,r)}$.
\end{definition}
\begin{remark}
A quadratic real space with signature $(p,q,0)$ is abbreviated by $\mathbb{R}^{(p,q)}$.
\end{remark}
\noindent
With the new basis $\lbrace e_1,\dots,e_n \rbrace$ the relations \eqref{EQ1} become
\begin{equation*}
e_ie_j+e_je_i=0,\quad i\neq j.
\end{equation*}
In the remainder of this paper we shall abbreviate the product of basis elements with lists
\begin{equation*}
e_{12\dots k}:=e_1e_2\dots e_k,\, \mbox{with $0\leq k\leq n$}.
\end{equation*}
The $2^n$ monomials 
\begin{equation*}
e_{i_1}e_{i_2}\dots e_{i_k},\quad 0\leq k\leq n
\end{equation*}
form the standard basis of the Clifford algebra. Furthermore, a Clifford algebra is the direct sum $\bigoplus\limits_{i=0}^n{\bigwedge}^i V$ of all exterior products ${\bigwedge}^i V$ of any grade $0\leq i\leq n$ where $e_{k_1}\dots e_{k_i},\,k_1<\dots<k_i$ form a basis of ${\bigwedge}^i V$. Therefore, a Clifford algebra is a graded algebra and its dimension is calculated by
\begin{equation*}
\dim \clifford_{(p,q,r)}=\sum\limits_{i=0}^n{\dim {\bigwedge}^i V}=\sum\limits_{i=0}^n\binom{n}{i}=2^n.
\end{equation*}
Moreover, the Clifford algebra $\clifford_{(p,q,r)}$ can be decomposed in an even and an odd part
\begin{align*}
\clifford_{(p,q,r)}&=\clifford_{(p,q,r)}^+ \oplus \clifford_{(p,q,r)}^-\\
&=\bigoplus\limits_{\substack{i=0\\i \text{ even}}}^n{{\bigwedge}^i V} \oplus \bigoplus\limits_{\substack{i=0\\i \text{ odd}}}^n{{\bigwedge}^i V}.
\end{align*}
The even part $\clifford_{(p,q,r)}^+$ is a subalgebra, because the product of two even graded monomials must be even graded since the generators cancel only in pairs. Elements contained in ${\bigwedge}^i V$ are called \emph{$i$-vectors} and the $\mathbb{R}$-linear combination of $i$-vectors $i\leq n$ is called a \emph{multi-vector}. A multi-vector $\mathfrak{A}$ is called homogeneous if $\left[ \mathfrak{A}\right]_i=\mathfrak{A},\, i\leq n$ where $\left[\cdot\right]_m\ m\in\mathbb{N}$ denotes the grade-$m$ part of the multi-vector $\mathfrak{A}$. The product of invertible vectors is called a \emph{versor} and the product of vectors where at least one vector is not invertible is called a \emph{null versor}.
\begin{definition} \label{centerclifford}
The center of a ring $\mathcal{R}$ is the set of all elements that commute with all other elements
\[\mathcal{C}(\mathcal{R}):=\left\lbrace c\in\mathcal{R}\mid cx=xc \mbox{ for all }a\in\mathcal{R}\right\rbrace.\]
\end{definition}
\noindent We are interested in the center of a Clifford algebra, see \cite[p.95]{garling:cliffordalgebras}. The center of a Clifford algebra $\clifford_{(p,q,r)}$ is
\begin{enumerate}
\item[(1)] $\mathcal{C}(\clifford_{(p,q,r)})=\lbrace\alpha e_0+\beta e_{12\ldots n}\mid\,\alpha,\beta\in\mathds{R}\rbrace$ if $n$ is odd,
\item[(2)] $\mathcal{C}(\clifford_{(p,q,r)})=\lbrace\alpha e_0\mid\,\alpha\in \mathds{R}\rbrace$ if $n$ is even.
\end{enumerate}
For the even part the center is
\begin{enumerate}
\item[(3)] $\mathcal{C}(\clifford_{(p,q,r)}^+)=\lbrace\alpha e_0\mid\,\alpha\in\mathds{R}\rbrace$ if $n$ is odd,
\item[(4)] $\mathcal{C}(\clifford_{(p,q,r)}^+)=\lbrace\alpha e_0+\beta e_{12\ldots n}\mid\,\alpha,\beta\in\mathds{R}\rbrace$ if $n$ is even.
\end{enumerate}
\subsection{Clifford Algebra Automorphisms}
\noindent For our purposes two automorphisms that exist on each Clifford algebra are interesting. The \emph{conjugation} is an \emph{anti-involution} denoted by an asterisk, see \cite{selig:geometricfundamentalsofrobotics}. Its effect on generators is given by $e_i^\ast=-e_i$. There is no effect on scalars. Extending the conjugation by using linearity yields
\begin{equation*}
(e_{i_1}e_{i_2}\dots e_{i_k})^\ast=(-1)^k e_{i_k}\dots e_{i_2}e_{i_1}
\end{equation*}
with $ 0\leq i_1<i_2<\dots <i_k\leq n$. The geometric product of a vector $\mathfrak{v}=\sum\limits_{i=1}^{n}{x_ie_i}\in{\bigwedge}^1 V$ with its conjugate results in
\begin{align*}
\mathfrak{v}\mathfrak{v}^\ast&=-x_1^2-x_2^2-\dots-x_p^2+x_{p+1}^2+\dots x_{p+q}^2\\
&=-b(v,v),
\end{align*}
where $v=(x_1,\dots x_n)^\mathrm{T}\in\mathbb{R}^n$.
\begin{definition}
The inverse element of a versor $\mathfrak{v}\in\clifford_{(p,q,r)}$ is defined by
\begin{equation*}
\mathfrak{v}^{-1}:=\frac{\mathfrak{v}^\ast}{N(\mathfrak{v})},
\end{equation*}
with $N(\mathfrak{v}):=\mathfrak{v}\mathfrak{v}^\ast$.
\end{definition}
\noindent
The map $N:\clifford_{(p,q,r)}\to\clifford_{(p,q,r)}$ is called the \emph{norm} of the Clifford algebra. For general multi-vectors $\mathfrak{M}\in\clifford_{(p,q,r)}$ inverse elements exist and are defined through the relation $\mathfrak{M}\mathfrak{M}^{-1}=\mathfrak{M}^{-1}\mathfrak{M}=1$, but the determination is more difficult and can be found in \cite{Fontijne:EfficientImplementationofGeometricAlgebra}. Note that in general not every element is invertible.
The other automorphism we are dealing with is the \emph{main involution}. It is denoted by $\alpha$ and defined by
\begin{equation*}
\alpha(e_{i_1}e_{i_2}\dots e_{i_k})=(-1)^k e_{i_1}e_{i_2}\dots e_{i_k}
\end{equation*}
for $ 0\leq i_1<i_2<\dots <i_k\leq n$. The main involution has no effect on the even subalgebra and it commutes with the conjugation, {\it i.e.}, $\alpha(\mathfrak{M}^\ast)=\alpha(\mathfrak{\mathfrak{M}})^\ast$ for arbitrary $\mathfrak{\mathfrak{M}}\in\clifford_{(p,q,r)}$.
\subsection{Clifford Algebra Products}
\noindent On vectors $\mathfrak{a},\mathfrak{b}\in{\bigwedge}^1 V$ we can write the inner product in terms of the geometric product
\begin{equation}\label{EQ10}
\mathfrak{a}\cdot \mathfrak{b}:=\frac{1}{2}(\mathfrak{a}\mathfrak{b}+\mathfrak{b}\mathfrak{a}).
\end{equation}
A generalization of the inner product to homogeneous multi-vectors can be found in \cite{hestenes:cliffordalgebra}.
For $\mathfrak{A}\in{\bigwedge}^k V,\ \mathfrak{B}\in{\bigwedge}^l V$ the generalized inner product is defined by
\begin{equation*}
\mathfrak{A}\cdot \mathfrak{B}:=\left[ \mathfrak{A}\mathfrak{B} \right]_{\vert k-l\vert}.
\end{equation*}
There is another product on vectors, {\it i.e.}, the outer (or exterior) product
\begin{equation}\label{EQ12}
\mathfrak{a}\wedge \mathfrak{b}:=\frac{1}{2}(\mathfrak{a}\mathfrak{b}-\mathfrak{b}\mathfrak{a}).
\end{equation}
This product can also be generalized to homogeneous multi-vectors, see again \cite{hestenes:cliffordalgebra}. For $\mathfrak{A}\in{\bigwedge}^k V,\ \mathfrak{B}\in{\bigwedge}^l V$ the generalized outer product is defined by
\begin{equation*}
\mathfrak{A}\wedge \mathfrak{B}:=\left[ \mathfrak{A}\mathfrak{B} \right]_{\vert k+l\vert}.
\end{equation*}
From equation \eqref{EQ10} and \eqref{EQ12} it follows that for vectors the geometric product can be written as the sum of the inner and the outer product
\begin{equation*}
\mathfrak{ab}=\mathfrak{a}\cdot \mathfrak{b}+\mathfrak{a}\wedge \mathfrak{b}.
\end{equation*}
More general this can be defined for multivectors with the commutator and the anti-commutator product, see \cite{perwass}. For treating geometric entities within this algebra context the definition of a $k$-blade, the \emph{inner product null space}, and its dual the \emph{outer product null space} is needed.
\begin{definition}
A $k$-blade is the $k$-fold exterior product of vectors $\mathfrak{v}\in \bigwedge^1 V$.
Therefore, a $k$-blade can be written as
\[\mathfrak{A}=\mathfrak{a}_1\wedge \mathfrak{a}_2 \wedge\dots \mathfrak{a}_k.\]
If the $k$-blade squares to zero, it is called a null $k$-blade.
The \emph{inner product null space} (IPNS) of a blade $\mathfrak{A}\in{\bigwedge}^k V$, cf. \cite{perwass}, is defined by
\begin{equation*}
\mathds{NI}(\mathfrak{A}):=\left\lbrace \mathfrak{v}\in{\bigwedge}^1 V\mid\mathfrak{v}\cdot \mathfrak{A}=0 \right\rbrace.
\end{equation*}
Moreover, the \emph{outer product null space} (OPNS) of a blade $\mathfrak{A}\in{\bigwedge}^k V$ is defined by
\begin{equation*}
\mathds{NO}(\mathfrak{A}):=\left\lbrace \mathfrak{v}\in{\bigwedge}^1 V\mid\mathfrak{v}\wedge \mathfrak{A}=0 \right\rbrace.
\end{equation*}
\end{definition}
\begin{remark}
The same set can be described with inner product or outer product null spaces. In the non-degenrate case the change between both representations is achieved with the \emph{pseudoscalar} $\mathfrak{J}=e_1\dots e_n$
\[\mathds{NO}(\mathfrak{A})=\mathds{NI}(\mathfrak{AJ}).\]
\end{remark}
\noindent Later on we work in a homogeneous Clifford algebra model, and therefore, it does not matter from which side we multiply with the pseudoscalar, since the inner product respectively outer product null space is not affected.
Moreover, multiplication with the pseudoscalar corresponds to the application of the polarity corresponding to the measure quadric. In the degenerate case the \emph{Poincar\'{e}-identity} has to be used, compare to \cite{gunn:kinematics}.
\subsection{Pin and Spin Groups}
\noindent With respect to the geometric product the units, {\it i.e.}, the invertible elements of a Clifford algebra denoted by $\clifford_{(p,q,r)}^\times$ form a group.
\begin{definition}
The \emph{Clifford group} is defined by
\begin{align*}
\Gamma(\clifford_{(p,q,r)})\!:=\!&\left\lbrace \mathfrak{g}\in\clifford_{(p,q,r)}^\times\mid \alpha(\mathfrak{g})\mathfrak{v} \mathfrak{g}^{-1}\in {\bigwedge}^1 V \right.  \\
 &\left. \mbox{ for all } \mathfrak{v}\in{\bigwedge}^1 V \right\rbrace.
\end{align*}
\end{definition}
\noindent
A proof that $\Gamma(\clifford_{(p,q,r)})$ is indeed a group with respect to the geometric product can be found in \cite{gallier:cliffordalgebrascliffordgroups}.
We define two important subgroups of the Clifford group.
\begin{definition}
The \emph{Pin group} is the subgroup of the Clifford group with $N(\mathfrak{g})=\pm 1$.
\begin{align*}
\mbox{Pin}_{(p,q,r)}\!:=\!&\left\lbrace \mathfrak{g}\in\clifford_{(p,q,r)}\mid \mathfrak{gg}^\ast\!=\!\pm 1\mbox{ and }\right.\\
 &\left.\alpha(\mathfrak{g})\mathfrak{v}\mathfrak{g}^\ast\in {\bigwedge}^1 V \mbox{ for all } \mathfrak{v}\in{\bigwedge}^1 V \right\rbrace.
\end{align*}
Furthermore, we define the \emph{Spin group} by $Pin_{(p,q,r)}\cap \clifford_{(p,q,r)}^+$
\begin{align*}
\mbox{Spin}_{(p,q,r)}\!:=\!&\left\lbrace \mathfrak{g}\in\clifford_{(p,q,r)}^+\mid \mathfrak{gg}^\ast\!=\!\pm 1\mbox{ and }\right.\\
&\left. \alpha(\mathfrak{g})\mathfrak{v}\mathfrak{g}^\ast\in {\bigwedge}^1 V \mbox{ for all } \mathfrak{v}\in{\bigwedge}^1 V \right\rbrace.
\end{align*}
\end{definition}
\begin{remark}
For non-degenerated Clifford algebras the Pin group is a double cover of the orthogonal group of the quadratic space $(V,b)$. Moreover, the Spin group is a double cover of the special orthogonal group of $(V,b)$, see \cite{gallier:cliffordalgebrascliffordgroups}.
\end{remark}

\section{The homogeneous Clifford Algebra model corresponding to Klein's Quadric}
To construct a homogeneous Clifford algebra model we use a vector space as model for the projective space which is the base space of the model. Thus, we use $\mathbb{R}^6$ as model space for $\mathbb{P}^5(\mathbb{R})$ together with the quadric form $b$ given by the corresponding symmetric matrix
\[\mathrm{B}=\begin{pmatrix}
\mathrm{O} & \mathrm{I}\\\mathrm{I} & \mathrm{O}
\end{pmatrix}\]
where $\mathrm{O}$ denotes the $3\times 3$ zero matrix and $\mathrm{I}$ the $3\times 3$ identity matrix. The $\bigwedge^1 V$ subspace is identified with the vector space, and therefore, with the projective space. The square of a general vector $\mathfrak{v}=\sum_ {i=1}^6{x_ie_i}$ results in
\[\mathfrak{vv}=x_1x_4+x_2x_5+x_3x_6.\]
This is the equation of Klein's quadric and a vector corresponds to a point on Klein's quadric if it is a null vector, {\it i.e.}, a vector squaring to zero. Since the Clifford algebra allows the description of subspaces with the use of blades of higher grade, we can express linear line manifolds within the homogeneous Clifford algebra model, see \cite{klawitter_lingeo, klawitter_diss}.
\subsection{Linear Line Manifolds}
Clifford algebras carry the subspace structure of Grassmann algebras together with the metric properties derived with the corresponding Cayley-Klein geometry. Therefore, we construct blades of grade $k$, $k\leq 5$ and examine the geometric inner product and geometric outer product null spaces.
\paragraph{Two-blades}
The exterior product of two vectors corresponds to a two-blade that describes a line $P_1^1$ in $\mathbb{P}^5(\mathbb{R})$. We use the outer product null space to determine the set of points on $P_1^1$ contained in $M_2^4$, and therefore, the set of lines in $\mathbb{P}^3(\mathbb{R})$ corresponding to the two-blade. If the whole line $P_1^1$ is contained in $M_2^4$ the corresponding two-blade is a null two-blade. Thus, the outer product null space corresponds to a pencil of lines in $\mathbb{P}^3(\mathbb{R})$ or to two lines in $\mathbb{P}^3(\mathbb{R})$ that may be identical. The inner product null space of a two-blade corresponds to a three-space intersection with Klein's quadric, {\it i.e.}, a linear line congruence in $\mathbb{P}^3(\mathbb{R})$.
\paragraph{Three-blades}
Three-blades can be generated as the exterior product of three vectors and describe two-spaces in $\mathbb{P}^5(\mathbb{R})$. The outer product null space of a three-blade corresponds to the two-space intersection of the two-space $P_1^2$ described by the three-blade with $M_2^4$. If the three-blade is a null three-blade $P_1^2$ is completely contained in Klein's quadric and the outer product null space corresponds to a bundle or a field of lines in $\mathbb{P}^3(\mathbb{R})$. The outer product null space of a non-null three-blade corresponds to a conic section of $M_2^4$ that defines a regulus in $\mathbb{P}^3(\mathbb{R})$. In a dual way the inner product null space of a three-blade is the intersection of the image of $P_1^2$ under the polarity corresponding to the meassure quadric with Klein's quadric.
\paragraph{Four-blades}
Outer product null spaces of four-blades correspond to three-space intersections with Klein's quadric and the corresponding set of lines in $\mathbb{P}^3(\mathbb{R})$ is a linear congruence of lines. The dual of a three-space in $\mathbb{P}^5(\mathbb{R})$ is a line. Thus, the inner product null space of a four-blade is the intersection of this line with $M_2^4$.
\paragraph{Five-blades}
Five-blades describe four-spaces in Klein's model and the intersection of a four-space with Klein's quadric results in a linear congruence of lines. The outer product null space of a five-blade corresponds to a linear congruence of lines. Moreover, the inner produt null space of a five-blade yields a vector.
\subsection{Clifford Group Action}\label{sectionCliffordgroup}
The null vectors correspond to points on Klein's quadric. Moreover, the sandwich action of an arbitrary versor on a null vector results a null vector again. Thus Klein's quadric is mapped to itself. To examine the action of the whole Clifford group, we first investigate the action of non-null vectors on null vectors.
\begin{remark}\label{remarksandwich}
Since we are working in a homogeneous Clifford algebra model multiplication with a real factor does not change the geometric meaning of an algebra element. Therefore, we use the conjugate element instead of the inverse element. Hence, the sandwich operator that we use is given by:
\[\alpha(\mathfrak{a})\mathfrak{va}^\ast=\mathfrak{aa}^\ast(\alpha(\mathfrak{a})\mathfrak{va}^{-1}).\]
This operator does not involve an inverse, and therefore, it can also be applied if the element $\mathfrak{a}$ is not invertible.
\end{remark}
\paragraph{Sandwich action of vectors}
Let $\mathfrak{a}=\sum_{i=1}^6{a_ie_i}$ and $\mathfrak{v}=\sum_{i=1}^6{x_ie_i}$
%\begin{align*}
%\mathfrak{a}&=a_1e_1\!+\!a_2e_2\!+\!a_3e_3\!+\!a_4e_4\!+\!a_5e_5\!+\!a_6e_6,\\
%\mathfrak{v}&=x_1e_1\!+\!x_2e_2\!+\!x_3e_3\!+\!x_4e_4\!+\!x_5e_5\!+\!x_6e_6
%\end{align*}
be two vectors with $\mathfrak{aa}\neq0$ and $\mathfrak{vv}=0$.
%\begin{align*}
%\mathfrak{aa}&=a_1a_4+a_2a_5+a_3a_6\neq 0,\\
%\mathfrak{vv}&=x_1x_4+x_2x_5+x_3x_6= 0.
%\end{align*}
The action of the sandwich operator $\alpha(\mathfrak{a})\mathfrak{va}^\ast$ on the $\bigwedge^1 V$ subspace can be expressed as product of a matrix with a vector. The matrix has the form:
{\arraycolsep4pt
\begin{equation}
\mathrm{M}=\begin{pmatrix}
  k_1 & a_1a_5   & a_1a_6 & a_1a_1 & a_1a_2 & a_1a_3\\
 a_2a_4 & k_2    & a_2a_6 & a_2a_1 & a_2a_2 & a_2a_3\\
 a_3a_4 & a_3a_5 & k_3    & a_3a_1 & a_3a_2 & a_3a_3 \\
 a_4a_4 & a_4a_5 & a_4a_6 & k_4    & a_4a_2 & a_4a_3\\
 a_5a_4 & a_5a_5 & a_5a_6 & a_5a_1 & k_5    & a_5a_3\\
 a_6a_4 & a_6a_5 & a_6a_6 & a_6a_1 & a_6a_2 &  k_6 
\end{pmatrix},
\end{equation}}
with
\begin{align*}
k_1&=-\!a_5a_2\!-\!a_6a_3,\quad &&k_2=-\!a_6a_3\!-\!a_4a_1,\\
k_3&=-\!a_4a_1\!-\!a_5a_2,\quad &&k_4=-\!a_5a_2\!-\!a_6a_3,\\
k_5&=-\!a_6a_3\!-\!a_4a_1,\quad &&k_6=-\!a_4a_1\!-\!a_5a_2.
\end{align*}
This involutoric automorphism of Klein's quadric corresponds to a null polarity acting on $\mathbb{P}^3(\mathbb{R})$, see \cite{klawitter_lingeo}. Furthermore, the Clifford group is generated by invertible  vectors, and therefore, all elements of the Clifford group correspond to compositions of null polarities in $\mathbb{P}^3(\mathbb{R})$. Moreover, this means that every element of $\clifford_{(3,3)}^+$ that is an even product of vectors corresponds to a collineation and every element of $\clifford_{(3,3)}^-$ that is an odd product of vectors corresponds to a correlation.

\paragraph{Collineations}
A general element $\mathfrak{g}\in\clifford_{(3,3)}^+$ corresponding to a collineation is given by
\begin{align*}
\mathfrak{g}&=g_{1}e_0\!+\!g_{2}e_{12}\!+\!g_{3}e_{13} \!+\!g_{4}e_{14}\!+\!g_{5}e_{15}\!+\!g_{6}e_{16}\\
&+\!g_{7}e_{23}\!+\!g_{8}e_{24}\!+\!g_{9}e_{25}\!+\!g_{10}e_{26}\!+g_{11}e_{34}\\
&+\!g_{12}e_{35}\!+\!g_{13}e_{36}\!+\!g_{14}e_{45}\!+\!g_{15}e_{46}\!+\!g_{16}e_{56}\\
&+\!g_{17}e_{1234}\!+\!g_{18}e_{1235}\!+\!g_{19}e_{1236}\!+\! g_{20}e_{1245}\\&
+g_{21}e_{1246}\!+\!g_{22}e_{1256} \!+\!g_{23}e_{1345}\!+\!g_{24}e_{1346}\\
&+g_{25}e_{1356}\!+\!g_{26}e_{1456}\!+\!g_{27}e_{2345}\!+\!g_{28}e_{2346}\\
&+g_{29}e_{2356}\!+\!g_{30}e_{2456}\!+\!g_{31}e_{3456}\!+\!g_{32}e_{123456}.
\end{align*}
We derive constraints to this element by $\alpha(\mathfrak{g})\mathfrak{vg}^\ast\in\bigwedge^1 V$ for all $\mathfrak{v}\in\bigwedge^1 V$. If we try to compute the corresponding collineation we have to distinguish the type of the collineation, {\it i.e.}, if it acts on points or on planes.
\paragraph{Action on Points}
First we assume $\mathfrak{g}\in\clifford_{(3,3)}^+$ corresponds to a collineation acting on points.
The automorphic collineation of Klein's quadric induced by $\mathfrak{g}$ can be transferred to a collineation that acts on $\mathbb{P}^3(\mathbb{R})$. Therefore, the element $g\in\clifford_{(3,3)}^+$ is applied to a null three-blade corresponding to a bundle of lines under the constraint that $\mathfrak{g}$ is a versor, see \cite{klawitter_lingeo}. This results in the matrix $(m_{i,j}), i,j=0,\dots,3$ for the collineation with
\begin{align*}
m_{00}&=g_{1}\!-\!g_{20}\!-\!g_{24}\!-\!g_{32}\!-\!g_{29}\!+\!g_{9}\!+\!g_{4}\!+\!g_{13},\\
m_{11}&=g_{24}\!-\!g_{9}\!+\!g_{20}\!-\!g_{13}\!-\!g_{32}\!+\!g_{1}\!+\!g_{4}\!-\!g_{29},\\
m_{22}&=g_{1}\!-\!g_{13}\!-\!g_{32}\!-\!g_{4}\!+\!g_{29}\!+\!g_{9}\!-\!g_{24}\!+\!g_{20},\\
m_{33}&=g_{24}\!+\!g_{13}\!+\!g_{29}\!+\!g_{1}\!-\!g_{4}\!-\!g_{9}\!-\!g_{20}\!-\!g_{32},
\end{align*}
and
\begin{align*}
m_{01}&=2(g_{7}+g_{17}),&&m_{02}=2(g_{18}\!-\!g_{3}),\\
m_{03}&=2(g_{19}\!+\!g_{2}),&&m_{10}=-\!2(g_{26}\!+\!g_{16}),\\
m_{12}&= 2(g_{5}\!+\!g_{25}), &&m_{13}=2(g_{6}\!-\!g_{22}),\\
m_{20}&= 2(g_{15}\!-\!g_{30}),&&m_{21}=2(g_{8}\!+\!g_{28}),\\
m_{23}&=2(g_{21}\!+\!2g_{10}),&&m_{30}=-\!2(g_{31}\!+\!g_{14}),\\
m_{31}&=2(g_{11}\!-\!g_{27}), &&m_{32}=2(g_{23}\!+\!g_{12}).
\end{align*}
The system of equations 
\begin{equation}\label{systemcollineation}
(m_{i,j})=(c_{i,j}),\ i,j=0,\dots,3
\end{equation}
can be used to compute an element of $\clifford_{(3,3)}^+$ corresponding to a collineation acting on points with matrix representation $\mathrm{C}=(c_{i,j}),\ i,j=0,\dots,3$. Therefore, constraint equations derived from $\alpha(\mathfrak{g})\mathfrak{vg}^\ast\in\bigwedge^1 V$ and $\mathfrak{gg}^\ast=\pm 1$ are used.\\

\paragraph{Action on Planes} Moreover, the same automorphic collineation of Klein's quadric may correspond to a collineation in $\mathbb{P}^3(\mathbb{R})$ that maps planes to planes. This action can also be represented as matrix vector product. Therefore, the general element $\mathfrak{g}\in\clifford_{(3,3)}^+$ is applied to a null three-blade corresponding to a field of lines. The action can be expressed as collineation with matrix $(m_{i,j}),\,i,j,=0,\dots,3$
\begin{align*}
m_{00}&=g_{32}\!-\!g_{20}\!-\!g_{13}\!-\!g_{29}\!-\!g_{9}\!-\!g_{24}\!+\!g_{1}\!-\!g_{4},\\
m_{11}&=g_{1}\!+\!g_{9}\!+\!g_{20}\!+\!g_{24}\!+\!g_{13}\!-\!g_{29}\!+\!g_{32}\!-\!g_{4},\\
m_{22}&=g_{20}\!+\!g_{29}\!+\!g_{4}\!+\!g_{13}\!-\!g_{9}\!-\!g_{24}\!+\!g_{1}\!+\!g_{32},\\
m_{33}&=g_{9}\!+\!g_{24}\!+\!g_{29}\!-\!g_{13}\!+\!g_{1}\!+\!g_{4}\!-\!g_{20}\!+\!g_{32},
\end{align*}
and
\begin{align*}
m_{01}&=2(g_{16}-g_{26}),&&m_{02}=-2(g_{15}+g_{30}),\\
m_{03}&=2(g_{14}-g_{31}),&&m_{10}=2(g_{17}-g_{7}),\\
m_{12}&=2(g_{28}-g_{8}), &&m_{13}=-2(g_{27}+g_{11}),\\
m_{20}&=2(g_{3}+g_{18}) ,&&m_{21}=2(g_{25}-g_{5}),\\
m_{23}&=2(g_{23}-g_{12}),&&m_{30}=2(g_{19}-g_{2}),\\
m_{31}&=-2(g_{22}+g_{6}), &&m_{32}=2(g_{21}-g_{10}).
\end{align*}
Depending on the action of the collineation, {\it i.e.}, if it acts on points or planes the elements of the Spin group can be transferred to their matrix representations and vice versa.
\paragraph{Correlations}
General elements of the odd part $\clifford_{(3,3)}^-$ have the form
\begin{align*}
\mathfrak{h} &= h_{1}e_{1}\!+\!h_{2}e_{2}\!+\!h_{3}e_{3}\!+\!h_{4}e_{4}\!+\!h_{5}e_{5}\!+\!h_{6}e_{6}\\
&+h_{7}e_{123}\!+\!h_{8}e_{124}\!+\!h_{9}e_{125}\!+\!h_{10}e_{126}\!+\!h_{11}e_{134}\!\\
&+h_{12}e_{135}\!+\!h_{13}e_{136}\!+\!h_{14}e_{145}\!+\!h_{15}e_{146}\!+\!h_{16}e_{156}\\
&+h_{17}e_{234}\!+\!h_{18}e_{235}\!+\!h_{19}e_{236}\!+\!h_{20}e_{245}\!+\!h_{21}e_{246}\\
&+h_{22}e_{256}\!+\!h_{23}e_{345}\!+\!h_{24}e_{346}\!+\!h_{25}e_{356}\!+\!h_{26}e_{456}\\
&+h_{27}e_{12345}\!+\!h_{28}e_{12346}\!+\!h_{29}e_{12356}\!+\!h_{30}e_{12456}\\
&+h_{31}e_{13456}\!+\!h_{32}e_{23456}.
\end{align*}
Their action on points or planes of the three-dimensional projective space $\mathbb{P}^3(\mathbb{R})$ can be described by the sandwich action on null three-blades corresponding to two-spaces contained entirely in $M_2^4$ that correspond to bundles of lines or fields of lines.
\paragraph{Action on Planes}
The action of a general element $\mathfrak{h}\in\clifford_{(3,3)}^-$ applied to a null three-blade corresponding to a two-space contained in Klein's quadric that is the image of a filed of lines under the Klein map can be expressed as $4\times 4$ matrix $(m_{i,j}), i,j=0,\dots,3$ that determines a correlation with
$m_{00}\!=\!-2h_{7},\ m_{11}\!=\!-2h_{16},\ m_{22}\!=\!2h_{21},\ m_{33}\!=\!-2h_{23},$
and
%\begin{align*}
%m_{01}&\!=\!h_{9}\!+\!h_{13}\!-\!h_{29}\!+\!h_{1},\!\! &&m_{02}\!=\!h_{2}\!-\!h_{8}\!+\!h_{28}\!+\!h_{19},\\
%m_{03}&\!=\!h_{3}\!-\!h_{27}\!-\!h_{18}\!-\!h_{11},\!\! &&m_{10}\!=\!h_{29}\!+\!h_{13}\!+\!h_{9}\!-\!h_{1},\\
%m_{12}&\!=\!h_{6}\!-\!h_{22}\!+\!h_{15}\!+\!h_{30},\!\! &&m_{13}\!=\!h_{31}\!-\!h_{25}\!-\!h_{5}\!-\!h_{14},\\
%m_{20}&\!=\!h_{19}\!-\!h_{8}\!-\!h_{2}\!-\!h_{28} ,\!\! &&m_{21}\!=\!h_{15}\!-\!h_{30}\!-\!h_{6}\!-\!h_{22},\\
%m_{23}&\!=\!h_{4}\!-\!h_{20}\!+\!h_{32}\!+\!h_{24},\!\! &&m_{30}\!=\!h_{27}\!-\!h_{11}\!-\!h_{18}\!-\!h_{3},\\
%m_{31}&\!=\!h_{5}\!-\!h_{31}\!-\!h_{25}\!-\!h_{14},\!\! &&m_{32}\!=\!h_{24}\!-\!h_{4}\!-\!h_{32}\!-\!h_{20}.\\
%\end{align*}
\begin{align*}
m_{01}&=h_{9}+h_{13}-h_{29}+h_{1},\\
m_{02}&=h_{2}-h_{8}+h_{28}+h_{19},\\
m_{03}&=h_{3}-h_{27}-h_{18}-h_{11},\\
m_{10}&=h_{29}+h_{13}+h_{9}-h_{1},\\
m_{12}&=h_{6}-h_{22}+h_{15}+h_{30},\\
m_{13}&=h_{31}-h_{25}-h_{5}-h_{14},\\
m_{20}&=h_{19}-h_{8}-h_{2}-h_{28},\\
m_{21}&=h_{15}-h_{30}-h_{6}-h_{22},\\
m_{23}&=h_{4}-h_{20}+h_{32}+h_{24},\\
m_{30}&=h_{27}-h_{11}-h_{18}-h_{3},\\
m_{31}&=h_{5}-h_{31}-h_{25}-h_{14},\\
m_{32}&=h_{24}-h_{4}-h_{32}-h_{20}.
\end{align*}
\paragraph{Action on Points} The same element $\mathfrak{h}\in\clifford_{(3,3)}^-$ can be interpreted as correlation in $\mathbb{P}^3(\mathbb{R})$ that acts on points and maps them to planes. This action can be described by a product of a matrix $m_{(i,j)},i,j=0,\dots,3$ with a vector. This matrix can be expressed with the coefficients of $\mathfrak{h}\in\clifford_{(3,3)}^-$ and has the entries $m_{00}\!=\!2h_{26},\ m_{11}\!=\!2h_{17},\ m_{22}\!=\!-2h_{12},\ m_{33}\!=\!2h_{10},$
and
%\begin{align*}
%m_{01}&\!=\!h_{32}\!-\!h_{4}\!-\!h_{20}\!-\!h_{24},\!\! &&m_{02}\!=\!h_{14}\!-\!h_{31}\!-\!h_{25}\!-\!h_{5},\\
%m_{03}&\!=\!h_{30}\!+\!h_{15}\!+\!h_{22}\!-\!h_{6},\!\! &&m_{10}\!=\!h_{4}\!-\!h_{32}\!-\!h_{24}\!-\!h_{20},\\
%m_{12}&\!=\!h_{18}\!-\!h_{27}\!-\!h_{3}\!-\!h_{11},\!\! &&m_{13}\!=\!h_{2}\!+\!h_{8}\!+\!h_{19}\!-\!h_{28},\\
%m_{20}&\!=\!h_{31}\!+\!h_{14}\!-\!h_{25}\!+\!h_{5},\!\! &&m_{21}\!=\!h_{3}\!-\!h_{11}\!+\!h_{27}\!+\!h_{18},\\
%m_{23}&\!=\!h_{9}\!-\!h_{13}\!-\!h_{1}\!-\!h_{29},\!\! &&m_{30}\!=\!h_{15}\!-\!h_{30}\!+\!h_{6}\!+\!h_{22},\\
%m_{31}&\!=\!h_{8}\!-\!h_{2}\!+\!h_{28}\!+\!h_{19},\!\! &&m_{32}\!=\!h_{1}\!-\!h_{13}\!+\!h_{29}\!+\!h_{9}.\\
%\end{align*}
\begin{align*}
m_{01}&=h_{32}-h_{4}-h_{20}-h_{24},\\
m_{02}&=h_{14}-h_{31}-h_{25}-h_{5},\\
m_{03}&=h_{30}+h_{15}+h_{22}-h_{6},\\
m_{10}&=h_{4}-h_{32}-h_{24}-h_{20},\\
m_{12}&=h_{18}-h_{27}-h_{3}-h_{11},\\
m_{13}&=h_{2}+h_{8}+h_{19}-h_{28},\\
m_{20}&=h_{31}+h_{14}-h_{25}+h_{5},\\
m_{21}&=h_{3}-h_{11}+h_{27}+h_{18},\\
m_{23}&=h_{9}-h_{13}-h_{1}-h_{29},\\
m_{30}&=h_{15}-h_{30}+h_{6}+h_{22},\\
m_{31}&=h_{8}-h_{2}+h_{28}+h_{19},\\
m_{32}&=h_{1}-h_{13}+h_{29}+h_{9}.
\end{align*}
Depending on the action of the correlation, {\it i.e.}, if it acts on points or planes the elements of the Pin group can be transferred to their matrix representations and vice versa.
\section{A Factorization Algorithm and its Application}
At this point we recall a theorem from \cite{klawitter_lingeo}.
\begin{theorem}\label{theoremKlein}
Every regular collineation or correlation can be expressed as the product of six null polarities at the most.
\end{theorem}
\noindent For a proof we refer to \cite{klawitter_lingeo} or \cite{klawitter_diss}. With the help of a modified version of a factorization algorithm introduced in \cite[p. 107]{perwass} we factorize arbitrary regular $4\times 4$ matrices into the product of six skew symmetric matrices at the most. Let $\left[  \mathfrak{g}\right]_k$ denote the grade-$k$ part of the versor $\mathfrak{g}$. For example a versor $\mathfrak{g}$ with maximal grade four can be written as $\mathfrak{g}=\left[  \mathfrak{g}\right]_0+\left[  \mathfrak{g}\right]_2+\left[  \mathfrak{g}\right]_4$.The part of maximal grade $\left[  \mathfrak{g}\right]_{\mathrm{max}}$ is always a blade, see \cite{perwass}. The non-null vector $\mathfrak{v}\in\bigwedge^1 V$ that is contained in
$\mathds{NO}(\left[  \mathfrak{g}\right]_{\mathrm{max}})$ satisfies
\begin{equation}\label{gradereduction}
\mathfrak{vg}_{\mathrm{max}}=\mathfrak{v}\cdot\mathfrak{g}_{\mathrm{max}}.
\end{equation}
Hence, $\mathfrak{g}'=\mathfrak{gv}^{-1}$ reduces the grade of the maximum blade of $\mathfrak{g}$ by one. With remark \ref{remarksandwich} this can be simplified to $\mathfrak{g}'=\mathfrak{gv}^\ast$ for homogeneous Clifford algebra models. Moreover, for vectors $\mathfrak{v}\in\bigwedge^1 V$ we have $\mathfrak{v}^\ast=-\mathfrak{v}$, and therefore, we can use  $\mathfrak{g}'=\mathfrak{gv}$ to reduce the degree of the initial versor $\mathfrak{g}$. The result of this product is again a versor. Repeated application of this process results in a set of vectors $\mathfrak{v}_1,\dots,\mathfrak{v}_{\mathrm{max}}$ whose geometric product results in the versor $\mathfrak{g}$ except for a real factor. Note that this procedure may fail for null versors. We show this algorithm in an example
\begin{example}\label{ex1}
Let $\mathrm{K}\in \mathrm{PGL}(\mathbb{P}^3(\mathbb{R}))$ be a collineation of $\mathbb{P}^3(\mathbb{R})$ that acts on points
\[K=\begin{pmatrix}
1 & 0 & 3 & 0\\1 & 1 & 0 & 1\\1 & 2 & 1 & 0\\1 & 1 & 2 & 1
\end{pmatrix}.\]
To get a versor $\mathfrak{g}\in\clifford_{(3,3)}^+$ corresponding to this collineation we have to solve the system \eqref{systemcollineation}.
\begin{align*}
&2(g_{7}\!+\!g_{17})\!=\!0, & 2(g_{18}\!-\!g_{3})\!=\!3,\\
&2(g_{19}\!+\!g_{2})\!=\!0,&-2(g_{26}\!+\!g_{16})\!=\!1,\\
&2(g_{5}\!+\!g_{25})\!=\!0, & 2(g_{6}\!-\!g_{22})\!=\!1,\\
&2(g_{15}\!-\!g_{30})\!=\!1, &2(g_{8}\!+\!g_{28})\!=\!2,\\
&2(g_{21}\!+\!2g_{10})\!=\!0, &-2(g_{31}\!+\!g_{14})\!=\!1,\\
& 2(g_{11}\!-\!g_{27})\!=\!1, &2(g_{23}\!+\!g_{12})\!=\!2,
\end{align*} 
%\begin{align*}
%&0\!=\!2(g_{7}\!+\!g_{17}),   &&3\!=\! 2(g_{18}\!-\!g_{3}),\\
%&0\!=\!2(g_{19}\!+\!g_{2}),   &&1\!=\!-2(g_{26}\!+\!g_{16}),\\
%&0\!=\!2(g_{5}\!+\!g_{25}),   &&1\!=\! 2(g_{6}\!-\!g_{22}),\\
%&1\!=\!2(g_{15}\!-\!g_{30}),  &&2\!=\!2(g_{8}\!+\!g_{28}),\\
%&0\!=\!2(g_{21}\!+\!2g_{10}), &&1\!=\!-2(g_{31}\!+\!g_{14}),\\
%&1\!=\!2(g_{11}\!-\!g_{27}), &&2\!=\!2(g_{23}\!+\!g_{12}),
%\end{align*} 
\begin{align*}
&g_{1}\!-\!g_{20}\!-\!g_{24}\!-\!g_{32}\!-\!g_{29}\!+\!g_{9}\!+\!g_{4}\!+\!g_{13}\!=\!1,\\
&g_{24}\!-\!g_{9}\!+\!g_{20}\!-\!g_{13}\!-\!g_{32}\!+\!g_{1}\!+\!g_{4}\!-\!g_{29}\!=\!1,\\
&g_{1}\!-\!g_{13}\!-\!g_{32}\!-\!g_{4}\!+\!g_{29}\!+\!g_{9}\!-\!g_{24}\!+\!g_{20}\!=\!1,\\
&g_{24}\!+\!g_{13}\!+\!g_{29}\!+\!g_{1}\!-\!g_{4}\!-\!g_{9}\!-\!g_{20}\!-\!g_{32}\!=\!1.
\end{align*}
There are two possibilities to guarantee that the resulting versor is in the Spin group, {\it i.e.}, $\mathfrak{gg}^\ast=1$ or $\mathfrak{gg}^\ast=-1$. We compute both solutions and start with the constraint equations implied by eq. $\alpha(\mathfrak{g})\mathfrak{vg}^\ast\in\bigwedge^1 V$ for all $\mathfrak{v}\in\bigwedge^1 V$ and $\mathfrak{gg}^\ast=1$.
The corresponding Spin group element has the form:
\begin{align*}
\mathfrak{g_+}&={1\over 8\sqrt{2}} \big(7 e_{0}\! +\!6 e_{12}\!-\!6 e_{13}\! +\! e_{14}\!- \!2 e_{15}\!- \!6 e_{23}\\
&+6 e_{24}\! - \!e_{25}\! -\!2 e_{26}\! +\!2 e_{34}\!+\!6 e_{35}\!-\!5 e_{36}\! -\!4 e_{45}\\
&+2 e_{46}\! +\!6 e_{1234}\! -\! 4 e_{56}\! +\!6 e_{1235}\! -\!6 e_{1236}\\
&-5 e_{1245}\!+\!2 e_{1246}\!-\!4 e_{1256}\! +\!2 e_{1345}\!-\!e_{1346}\\
&+2 e_{1356}\!-\!2 e_{2345}\!+\!2 e_{2346} \!+\! e_{2356}\!-\!2 e_{2456} \\
&-e_{123456}\big).
\end{align*}
If we demand that $\mathfrak{gg}^\ast=-1$ the resulting Spin group element is computed as
\begin{align*}
\mathfrak{g}_-&={1\over 8\sqrt{2}}\big(e_{0}\!-\!6 e_{12}\!-\!6 e_{13}\!-\!e_{14}\!+\!2 e_{15}\!+\!4 e_{16}\\
&+6 e_{23}\!+\!2 e_{24}\!+\!e_{25}\!+\!2 e_{26}\!+\!2 e_{34}\!+\!2 e_{35}\!+\!5 e_{36}\\
&+2 e_{46}\!-\!6 e_{1234}\!+\!6 e_{1235}\!+\!6 e_{1236}\!+\!5 e_{1245}\\
&-2 e_{1246}\!+\!6 e_{1345}\!+\!e_{1346}\!-\!2 e_{1356}\!-\!4 e_{1456}\\
&-2 e_{2345}\!+\!6 e_{2346}\!-\!e_{2356}\!-\!2 e_{2456}\!-\!4 e_{3456}\\
&-7 e_{123456}\big).
\end{align*}
Both elements $\mathfrak{g}_+$ and $\mathfrak{g}_-$ correspond to the same collineation whose entries can be computed with the coefficients of $\mathfrak{g}_+$ and $\mathfrak{g}_-$ and $(m_{i,j}),i,j=0,\dots, 3$.
\end{example}
\noindent Depending on the matrix of the projective transformation, {\it i.e.}, if it describes a collineation or correlation that acts on points or planes the corresponding versor can be computed with the use of the corresponding systems derived in the last section.
\begin{remark}
The matrix 
\[K=\begin{pmatrix}
-1 & 0 & 3 & 0\\1 & 1 & 0 & 1\\1 & 2 & 1 & 0\\1 & 1 & 2 & 1
\end{pmatrix}\]
results in a versor with complex entries. Thus, we have to work over the complex numbers.
\end{remark}
\noindent Moreover, it can be verified that $\mathfrak{g}_+\corr \mathfrak{Jg}_-$, where $\corr$ means is equal up to a scalar factor. Multiplication with the pseudoscalar does not change the action of the versor since
\begin{align*}
(\mathfrak{Jg})\mathfrak{v}(\mathfrak{Jg})^\ast&=\alpha(\mathfrak{J})\alpha(\mathfrak{g})\mathfrak{v}\mathfrak{g}^\ast\mathfrak{J}^\ast\\
&=\mathfrak{J}\mathfrak{J}^\ast\alpha(\mathfrak{g})\mathfrak{v}\mathfrak{g}^\ast\\
&\corr \alpha(\mathfrak{g})\mathfrak{v}\mathfrak{g}^\ast.
\end{align*}
\begin{remark}
Multiplication with an element of the center $\mathcal{C}(\clifford_{(3,3)})=a e_0 + b e_{123456},\ a,b\in\mathbb{C}$ has no effect on the sandwich action of a versor. Therefore, we have a group isomorphism
\begin{large}
\[\nicefrac{\mathrm{PGL(4,\mathbb{C})}}{\mathcal{C}(\mathrm{PGL(4,\mathbb{C})})}\to\nicefrac{\mathrm{Pin}(\clifford_{(3,3)})}{\mathcal{C}(\clifford_{(3,3)})}\]
\end{large}
\end{remark}
\noindent Now we factorize the element $\mathfrak{g}_+$ into vectors that correspond to null polarities. Due to the homogeneous setting we start with the numerator of $\mathfrak{g}_+$ of ex. \ref{ex1} denoted by $\mathfrak{g}\ (=8\sqrt{2}\mathfrak{g}_+)$.
\begin{example}
To factorize $\mathfrak{g}$ we first observe that $\left[ \mathfrak{g}\right]_\mathrm{max}=\left[ \mathfrak{g}\right]_6=-e_{123456}$. Multiplication with a scalar does not change the outer product null space of a $k$-blade, and therefore, it follows
\[\mathds{NO}(\left[ \mathfrak{g}\right]_\mathrm{max})=\mathds{NO}(e_{123456})=\lbrace\mathfrak{v}\in{\bigwedge}^1 V\rbrace.\]
We can now chose an arbitrary non-null vector $\mathfrak{v}_1\in\mathds{NO}(\left[ \mathfrak{g}\right]_\mathrm{max})$. Let $\mathfrak{v}_1=e_1\!+\!e_4$
be the first factor. The versor $\mathfrak{g}_1\!=\!\mathfrak{gv}_1$ has maximal grade five and the grade-$5$ part is given by
\[\left[ \mathfrak{g}_1\right]_5\! \corr \!e_{23456}\!-\!4e_{12345}\!+\!4e_{12346}\!-\!3e_{12456}\!+\!e_{13456}.\]
The outer product null space of $\left[ \mathfrak{g}_1\right]_\mathrm{5}$ is computed as
\begin{equation*}
\mathds{NO}(\left[ \mathfrak{g}_1\right]_5)=\lbrace\mathfrak{a}\in{\bigwedge}^1 V \mid \mathfrak{a \wedge \left[ \mathfrak{g}_1\right]_5}=0\rbrace.
\end{equation*}
This results in the set of all vectors $\mathfrak{a}=\sum_{i=1}^6{a_ie_i}$
satisfying $a_1\!-\!a_2\!-\!3a_3\!+\!4a_5\!+\!4a_6\!=\!0$.
We chose the non-null vector $\mathfrak{v}_2=4 e_2\!+\! e_5\in\mathds{NO}(\left[ \mathfrak{g}_1\right]_5)$
and compute the grade-four part of $\mathfrak{g}_2\!=\!\mathfrak{g}_1\mathfrak{v}_2$
\begin{align*}
\left[ \mathfrak{g}_2\right]_4&\corr e_{2356}\!-\!12e_{1234}\!+\!5e_{2345}\!-\!4e_{1236}\!+\!4e_{1235}\\
&+\!e_{1245}\!+\!e_{1345}\!-\!8e_{2346}\!+\!4e_{2456}\!-\!e_{3456}\\
&+8e_{1246}\!-\!3e_{1256}\!-\!4e_{1346}\!+\!e_{1356}\!+\!e_{1456}.
\end{align*}
Again we determine a vector contained in the outer product null space of $\left[ \mathfrak{g}_2\right]_4$. This null space is three-dimensional and spanned by the vectors
\begin{align*}
\mathfrak{b}_1&=e_1\!+\!5 e_2 \!+\! e_6, &&\mathfrak{b}_2=4e_1\!+\!8e_2\!+\!e_5\\
\mathfrak{b}_3&=e_1\!+\!e_2\!+\!e_4, &&\mathfrak{b}_4=e_1\!+\!4e_2\!-\!e_3.
\end{align*}
We chose the non-null vector $\mathfrak{b}_3$ and set $\mathfrak{v}_3=\mathfrak{b}_3$. In the next step we have to compute the outer product null space of $\left[ \mathfrak{g}_3\right]_3$ with $\mathfrak{g}_3=\mathfrak{g}_2\mathfrak{v}_3$ that is two-dimensional and spanned by:
\begin{align*}
\mathfrak{b}_1&=e_1\!+\!4e_2\!-\!3e_3,&&\mathfrak{b}_2=3e_1\!+\!7e_2\!-\!e_4\!+\!e_5,\\
\mathfrak{b}_3&=4e_2\!-\!e_4\!+\!e_6
\end{align*}
To keep the factors simple we set
\[\mathfrak{v}_4=\mathfrak{b}_1+\mathfrak{b}_3=-e_1+e_3-e_4+e_6.\]
The reduced versor $\mathfrak{g}_4=\mathfrak{g}_3\mathfrak{v}_4$ has maximal grade two and the outer product null space of the grade-two part is spanned by the two vectors:
\begin{align*}
\mathfrak{b}_1&=2e_1\!+\!3e_2\!+\!e_3\!-\!e_4\!+\!e_5,\\
\mathfrak{b}_2&=4e_2\!-\!e_4\!+\!e_6.
\end{align*}
Note that the vector $\mathfrak{b}_2$ is a null vector, and therefore, we chose $\mathfrak{v}_5=\mathfrak{b}_1$ as the next factor. The last factor $\mathfrak{v}_6$ is obtained by
\begin{align*}
\mathfrak{v}_6&=\mathfrak{g}_4\mathfrak{v}_5\\
&=8e_1\!+\!8e_2\!+\!4e_3\!-\!3e_4\!+\!4e_5\!-\!e_6.
\end{align*}
This results in the factorization
\[\mathfrak{v}_6\mathfrak{v}_5\mathfrak{v}_4\mathfrak{v}_3\mathfrak{v}_2\mathfrak{v}_1\corr\mathfrak{g}.\]
The next step is to carry over the factorization to projective transformations. Each vector corresponds to a null polarity. Therefore, we aim at a matrix product of the form
\[\mathrm{M}_6\mathrm{M}_5\mathrm{M}_4\mathrm{M}_3\mathrm{M}_2\mathrm{M}_1\corr\mathrm{K}.\]
Since we started with a collineation that maps points to points we have to ensure that the vectors are transferred to projective mappings of the right type. For example the product $\mathrm{M}_2\mathrm{M}_1$ must be the matrix corresponding to a collineation that maps points to points. Hence, $\mathrm{M}_1$ has to correspond to a null polarity acting on points and $\mathrm{M}_2$ to a null polarity acting on planes. Thus, the matrix product has to be understand as
\[\mathrm{M}_6\mathrm{M}_5^\ast\mathrm{M}_4\mathrm{M}_3^\ast\mathrm{M}_2\mathrm{M}_1^\ast\corr\mathrm{K},\]
where $\mathrm{M}_i^\ast,\ i=1,3,5$ is a null polarity acting on points and $\mathrm{M}_i,\ i=2,4,6$ is a null polarity acting on planes. We compute the matrices $\mathrm{M}_i^\ast$ and $\mathrm{M}_i$ with the matrix representations computed in section \ref{sectionCliffordgroup}. This results in the six null polarities:
{\arraycolsep3pt
\begin{align*}
\mathrm{M}_1^\ast&\!\!=\!\!\begin{small}\begin{pmatrix} 0 & \!-\!1 & 0 & 0 \\1 & 0 & 0 & 0 \\0 & 0 & 0 & \!-\!1 \\0 & 0 & 1 & 0 \end{pmatrix}\end{small}\!\!\!,
&&\mathrm{M}_2\!\!=\!\!\begin{small}\begin{pmatrix} 0 & 0 & 4 & 0 \\0 & 0 & 0 & \!-\!1 \\\!-4 & 0 & 0 & 0 \\0 & 1 & 0 & 0\end{pmatrix}\end{small}\!\!\!,\\
\mathrm{M}_3^\ast&\!\!=\!\!\begin{small}\begin{pmatrix} 0 & \!-\!1 & 0 & 0 \\1 & 0 & 0 & 1 \\0 & 0 & 0 & \!-\!1 \\0 & \!-\!1 & 1 & 0\end{pmatrix}\end{small}\!\!\!,
&&\mathrm{M}_4\!\!=\!\!\begin{small}\begin{pmatrix} 0 & \!-\!1 & 0 & 1 \\1 & 0 & 1 & 0 \\0 & \!-\!1 & 0 & \!-\!1 \\\!-\!1 & 0 & 1 & 0\end{pmatrix}\end{small}\!\!\!,\\
\mathrm{M}_5^\ast&\!\!=\!\!\begin{small}\begin{pmatrix} 0 & 1 & \!-\!1 & 0 \\\!-\!1 & 0 & \!-\!1 & 3 \\1 & 1 & 0 & \!-\!2 \\0 & \!-\!3 & 2 & 0\end{pmatrix}\end{small}\!\!\!,
&&\mathrm{M}_6\!\!=\!\!\begin{small}\begin{pmatrix} 0 & 8 & 8 & 4 \\\!-\!8 & 0 & \!-\!1 & \!-4 \\\!-\!8 & 1 & 0 & \!-\!3 \\\!-4 & 4 & 3 & 0\end{pmatrix}\end{small}\!\!\!.
\end{align*}}

\noindent With these six skew symmetric matrices it can be verified that
\[\mathrm{M}_6\mathrm{M}_5^\ast\mathrm{M}_4\mathrm{M}_3^\ast\mathrm{M}_2\mathrm{M}_1^\ast=-{1\over 4}\mathrm{K}.\]
Thus, we have found one possible factorization of the initial projective mapping into six null polarities.
\end{example}
\noindent It is worth to check the maximal grade of a versor that is solution of the other systems that were derived in section \ref{sectionCliffordgroup} if a factorization of a $4\times 4$ matrix is searched. The maximal grade of the versor determines the number of skew-symmetric matrices that are necessary to factorize the given matrix and for correlations the maximal possible grade is five.
\section{The homogeneous Clifford Algebra Model of Lie's Sphere Geometry}
%The construction of a homogeneous Clifford algebra model can be performed for any quadric. One important topic is the homogeneous Clifford algebra model corresponding to Lie's sphere geometry.
The same construction that we described for Klein's quadric can be applied to any quadric. As a second example that shall demonstrate the power of this calculus we examine Lie sphere geometry in a Clifford algebra context.
Lie sphere geometry is the geometry of oriented spheres. Especially, for the three-dimensional case the set of oriented spheres can be mapped  to a hyperquadric $L^4_1$ in five-dimensional projective space $\mathbb{P}^5(\mathbb{R})$. The construction goes back to S. \textsc{Lie} and was treated again by W. \textsc{Blaschke}, cf. \cite{Blaschke:VorlesungenuberDifferentialgeometrieundGeometrischeGrundlagenvonEinsteinsRelativitatstheorieBd.3}. A modern treatment of this topic can be found in \cite{Cecil:Liespheregeometry}. Moreover, the Lie construction can be achieved for arbitrary dimension.

\subsection{Lie's Quadric}\index{Lie's quadric}\label{subsectionLiequadrik}
A point model for the set of oriented hyperspheres, hyperplanes, and points (considered as spheres of radius $0$) of $\mathbb{R}^n$ is given by the projective hyperquadric
\[L_1^{n+1}:\,-x_0^2+x_1^2+\ldots+x_{n+1}^2-x_{n+2}^2=0.\]
%\FloatBarrier
For our purposes it is convenient that we restrict ourselves to the case of oriented spheres in three-dimensional Euclidean space. Nevertheless, we formulate the calculus for arbitrary dimensions. The quadric $L^{n+1}_1\subset\mathbb{P}^{n+2}(\mathbb{R})$ is of dimension $n\!+\!1$, degree $2$, and is called Lie's quadric. The maximal dimension of subspaces contained by $L_1^{n+1}$ is $1$, and therefore, there are no two-spaces contained entirely in $L_1^{n+1}$. Oriented hyperspheres, hyperplanes, and points are represented in Lie coordinates as shown in Table \ref{talbeliecoords}.
%{  \renewcommand{\arraystretch}{1.4}
%\[\begin{array}{c|c}\mbox{\textbf{Euclidean}}&\mbox{\textbf{Lie}}\\\hline
%\mbox{points: $u\in\mathbb{R}^n$}&\left( \frac{1+u\cdot u}{2},\frac{1-u\cdot u}{2},u,0  \right)^\mathrm{T} \mathbb{R}\\
%\infty &  \left( 1,-1,0,0\right)^\mathrm{T} \mathbb{R} \\
%\mbox{sphere: center $p$, signed radius $r$} &  \left( \frac{1+p\cdot p-r^2}{2},\frac{1-p\cdot p+r^2}{2},p,r \right)^\mathrm{T} \mathbb{R} \\
%\mbox{planes: $u\cdot N=h$, unit normal $N$} &  \left( h,-h,N,1\right)^\mathrm{T} \mathbb{R} 
%\end{array}\]
%}
It is not difficult to recover the Euclidean representation from Lie coordinates. If $x_0+x_1=0$ and if $x_{n+2}=0$ we have the point at infinity. If $x_{n+2}\neq 0$ we bring the point to the form $ \left( h,-h,N,1\right)^\mathrm{T}\mathbb{R} $ by dividing by $x_{n+2}$. If $x_0+x_1\neq 0$ and if $x_{n+2}=0$, we have a proper point. We obtain its normal form by dividing by $x_0+x_1$. The last case is if $x_{n+2}\neq 0$. In this case we have an oriented sphere. Again we get its normal form through division by $x_0+x_1$.\\

{\linespread{1.4}
\begin{table*}[hbt]
\begin{center}
\begin{tabular}{|c|c|}
\hline 
\textbf{Euclidean}&\textbf{Lie}\\
\hline
points: $u\in\mathbb{R}^n$ & $\left( \frac{1+u\cdot u}{2},\frac{1-u\cdot u}{2},u_1,\dots,u_n,0  \right)^\mathrm{T} \mathbb{R}$\\
$\infty$ &  $\left( 1,-1,0,\dots,0,0\right)^\mathrm{T} \mathbb{R}$ \\
sphere: center $p\in\mathbb{R}^n$, signed radius $r$ &  $\left( \frac{1+p\cdot p-r^2}{2},\frac{1-p\cdot p+r^2}{2},p_1,\dots,p_n,r \right)^\mathrm{T} \mathbb{R}$ \\
planes: $u\cdot N=h$, unit normal $N\in\mathbb{R}^n$ &  $\left( h,-h,N_1,\dots,N_n,1\right)^\mathrm{T} \mathbb{R} $\\
\hline 
\end{tabular} 
\caption{Correspondence between Euclidean entities and Lie-coordinates.}
\label{talbeliecoords}
\end{center}
\end{table*}}
\noindent The fundamental invariant of Lie sphere geometry is the oriented contact of spheres. It is not difficult to show that two spheres are in oriented contact if, and only if, their Lie coordinates $s_1,s_2\in L_1^{n+1}$ satisfy $\ell(s_1,s_2) =0$, where $\ell(\cdot,\cdot)$ denotes the bilinear form corresponding to $L_1^{n+1}$.\\

\noindent Especially for $n=3$ the lines on $L^4_1$ correspond to so called parabolic pencils of spheres. These pencils consist of all oriented spheres with one common point of contact. Furthermore, each parabolic pencil contains exactly one point, {\it i.e.}, sphere of radius $0$. If this point sphere is not $\infty$ the pencil contains exactly one oriented hyperplane $\Sigma$.
\begin{remark}
Conics on Lie's quadric correspond to Dupin cyclides, that are the envelopes of two one-parameter families of spheres.
\end{remark}
\noindent The group of Lie transformations shows up as the group of projective automorphisms of $L^{n+1}_1$. This group is isomorphic to $\mathrm{O}(n+1,2)/\pm 1$, see \textsc{Cecil} \cite{Cecil:Liespheregeometry}. Since the Pin group of the Clifford algebra $\clifford_{(n+1,2,0)}$ is a double cover of $\mathrm{O}(n+1,2)$ we can use this group to describe Lie transformations.
%\FloatBarrier
\subsection{The homogeneous Clifford Algebra Model corresponding to Lie Sphere Geometry}\label{subsectionLieClifford}
In this section we discuss the Clifford algebra model for Lie Sphere Geometry in the three-dimensional case. Therefore, the projective space we are dealing with is a five-dimensional space $\mathbb{P}^5(\mathbb{R})$. The homogeneous Clifford algebra model is obtained with the six-dimensional real vector space $\mathbb{R}^6$ as a model for the projective image space together with the quadratic form of Lie's quadric
\[\mathrm{Q}=\begin{pmatrix}
-1 & 0 &0 & 0 &0 &0\\
 0 & 1 &0 & 0 &0 &0\\
 0 & 0 &1 & 0 &0 &0\\
 0 & 0 &0 & 1 &0 &0\\
 0 & 0 &0 & 0 &1 &0\\
 0 & 0 &0 & 0 &0 &-1\\
\end{pmatrix}.\]
This algebra has signature $(p,q,r)=(4,2,0)$ and is of dimension $2^6=64$. Again the advantage of the Clifford algebra lies in the common description of the application of Lie transformations. Arbitrary projective subspaces of $\mathbb{P}^5(\mathbb{R})$ are transformed by the sandwich operator. As an example we determine all Lie inversions that leave the point at infinity fixed, {\it i.e.}, the subgroup of \emph{Laguerre} transformations. The point at infinity has the form $\mathfrak{p}=e_1\!-\!e_2$, compare to Table \ref{talbeliecoords}. A general invertible vector is given by $\mathfrak{a}=\sum_{i=1}^6{a_ie_i}$
%\begin{align*}
%\mathfrak{a}&=a_1e_1\!+\!a_2e_2\!+\!a_3e_3\!+\!a_4e_4\!+\!a_5e_5\!+\!a_6e_6,\\
%&-\!a_1^2\!+\!a_2^2\!+\!a_3^2\!+\!a_4^2\!+\!a_5^2\!-\!a_6^2\neq 0.
%\end{align*}
with $\mathfrak{aa}\neq 0$. The application of the sandwich operator to $\mathfrak{p}$ results in
\begin{align*}
\alpha(\mathfrak{a})\mathfrak{pa}^\ast&=-2(a_1\!+\!a_2)a_3e_3\!-\!2(a_1\!+\!a_2)a_4e_4\\
-&(a_1^2\!+\!a_2^2\!+\!a_3^2\!+\!a_4^2\!+\!a_5^2\!-\!a_6^2\!+\!2a_2a_1)e_1\\
-&(a_1^2\!+\!a_2^2\!-\!a_3^2\!-\!a_4^2\!-\!a_5^2\!+\!a_6^2\!+\!2a_2a_1)e_2 \\
-&2(a_1\!+\!a_2)a_5e_5\!-\!2(a_1\!+\!a_2)a_6e_6.
\end{align*}
To guarantee that this entity represents the point at infinity, we first see that $a_1\!+\!a_2\!=\!0$. With this condition the coefficients of $e_3,e_4,e_5$, and $e_6$ vanish. Moreover, the sum of the coefficients $e_1$ and $e_2$ has to vanish. This results in
\[-\!2a_1^2\!-\!2a_2^2\!-\!4a_2a_1\!=\!-\!2(a_1\!+\!a_2)^2\!=\!0.\]
Therefore, the only condition to a Lie inversion that it represents a Laguerre transformation is given by $a_1\!+\!a_2\!=\!0$ and the subgroup of Laguerre transformations is generated by all vectors with $a_1\!+\!a_2\!=\!0$.
\begin{remark}
Analogue to Th. \ref{theoremKlein} we can formulate a similar theorem for Lie sphere geometry in arbitrary dimensions. Since the projective model space for $n$-dimensional Lie sphere geometry has dimension $n+2$, the vector space for the homogeneous Clifford algebra model has dimension $n+3$. That means the highest grade is equal to $n+3$, and therefore, every group element can be written as the composition of $n+3$ vectors at the most. Especially for the case $n=3$, we have similar results as for Klein's quadric. In this case six involutions are necessary to generate the whole group.
\end{remark}
\noindent Let us reformulate this remark as theorem.
\begin{theorem}
Every Lie transformation in $n$-dimensional space is the composition of $n\!+\!3$ involutions, that correspond to the sandwich action of vectors.
\end{theorem}

\section{Conclusion}
With the help of the homogeneous Clifford algebra model corresponding to line geometry we introduced null polarities as fundamental involutions. After we transferred linear line manifolds to this model, we derived the correspondence between regular projective transformations and versors. Every versor can be factorized into vectors that correspond to null polarities. Therefore, every $4\times 4$ matrix corresponding to a projective transformation can be expressed as the product of six skew symmetric matrices at the most. Moreover, we presented a homogeneous Clifford algebra model corresponding to Lie sphere geometry.
\section*{Acknowledgments}

This work was supported by the research project "Line Geometry for Lightweight Structures", funded by the DFG 
(German Research Foundation) as part of the SPP 1542.

%%% You can use BibTeX for formatting your references:
%\bibliography{icgg}

%%% If you prefer to format your references manually, do it as
%%% follows:
\setlength{\bibsep}{2\parskip}

% Every article should conclude with some basic information about the
% authors:
\section*{About the author}
\begin{enumerate}
\item Daniel Klawitter, Dipl.-Math., is a PhD student at the Institute of Geometry at Dresden University of Technology, Germany. His research interests are Clifford algebras, line geometry, computational and theoretical kinematics, differential geometry and computer aided geometric design. He can be reached by e-mail: daniel.klawitter@tu-dresden.de or through postal address: TU Dresden, Institute of Geometry, 01062 Dresden, Germany.
\end{enumerate}

\end{document}